\documentclass[11pt,a4paper]{article}
\usepackage[latin1]{inputenc}
\usepackage{graphicx}
\usepackage[mathscr]{euscript}
\usepackage{a4}
\usepackage{epsfig}
\usepackage{stmaryrd,amsfonts,amsmath,amssymb}
\usepackage{ifthen}
\usepackage{hyperref}

\newcommand{\origforall}{}    \let\origforall=\forall
\renewcommand{\forall}{\,\origforall\,}
\newcommand{\origexists}{}    \let\origexists=\exists
\renewcommand{\exists}{\,\origexists\,}
\newcommand{\origin}{}        \let\origin=\in
\renewcommand{\in}{{\,\origin\,}}
\newcommand{\origwedge}{}     \let\origwedge=\wedge
\renewcommand{\wedge}{{\quad\origwedge\quad}}
\newcommand{\origvee}{}       \let\origvee=\vee
\renewcommand{\vee}{{\quad\origvee\quad}}
\newcommand{\origfrac}{}      \let\origfrac=\frac
\renewcommand{\frac}[2]{{\;\origfrac{#1}{#2}\;}}
    
\renewcommand{\subset}{\subseteq}
    
\renewcommand{\supset}{\supseteq}

\newcommand{\beq}{\begin{eqnarray*}}
\newcommand{\eeq}{\end{eqnarray*}}

\renewcommand{\d}{\n{d}}

\newcommand{\R}{\mathbb{R}}

\newcommand{\N}{\mathbb{N}}
\newcommand{\Z}{\mathbb{Z}}

\newcommand{\n}{\textnormal}

\DeclareMathOperator{\Lip}{{\n{Lip}}}
\newcommand{\smax}{\origvee}
\newcommand{\smin}{\origwedge}
\newcommand{\bigsup}{\bigvee}
\newcommand{\biginf}{\bigwedge}

\DeclareMathOperator{\lcm}{\n{lcm}}
\DeclareMathOperator{\PG}{\n{PG}}
\DeclareMathOperator{\mli}{\n{cmli}}
\DeclareMathOperator{\hyp}{\n{hyp}}
\DeclareMathOperator{\Pot}{\wp}
\DeclareMathOperator{\im}{\n{im}}
\newcommand{\my}{\mu}
\newcommand{\sd}{\vartriangle}
\newcommand{\imetric}{me\-tric~}  
\newcommand{\down}[1]{\Downarrow\! #1}

\newcommand{\esquare}{\hfill\ensuremath{\square}}

\newcounter{claim}
\newcounter{titleclaim}
\newenvironment{Proof}{{\textbf{Proof}$\;$}}{\esquare\newline}
\newenvironment{Non-Proof}{{\textbf{Non-Proof}$\;$}}{\hfill\ensuremath{\times}\newline}

\newenvironment{Question}{{\textbf{Question}$\;$}\begin{itshape}}{\end{itshape}}

\newtheorem{LemmaBase}[claim]{Lemma}
\newenvironment{Lemma}[1][-1]{\ifthenelse{\equal{#1}{-1}}
    {\begin{LemmaBase}}{\begin{LemmaBase}[#1]}
    \dotfill \, {\em \bf \arabic{claim}} \quad\\{} }
    {\addtocounter{titleclaim}{1}\end{LemmaBase}}

\newtheorem{CorollaryBase}[claim]{Corollary}
\newenvironment{Corollary}[1][-1]{\ifthenelse{\equal{#1}{-1}}
    {\begin{CorollaryBase}}{\begin{CorollaryBase}[#1]}
    \dotfill \, {\em \bf \arabic{claim}} \quad\\{} }
    {\addtocounter{titleclaim}{1}\end{CorollaryBase}}

\newtheorem{PropositionBase}[claim]{Proposition}
\newenvironment{Proposition}[1][-1]{\ifthenelse{\equal{#1}{-1}}
    {\begin{PropositionBase}}{\begin{PropositionBase}[#1]}
    \dotfill \, {\em \bf \arabic{claim}} \quad\\{} }
    {\addtocounter{titleclaim}{1}\end{PropositionBase}}

\newtheorem{ConjectureBase}[claim]{Conjecture}

\newtheorem{ExampleBase}[claim]{Example}
\newenvironment{Example}[1][-1]{\ifthenelse{\equal{#1}{-1}}
    {\begin{ExampleBase}}{\begin{ExampleBase}[#1]}
    \dotfill \, {\em \bf \arabic{claim}} \quad\\{} }
    {\addtocounter{titleclaim}{1}\end{ExampleBase}}

\newtheorem{DefinitionBase}[claim]{Definition}
\newenvironment{Definition}[1][-1]{\ifthenelse{\equal{#1}{-1}}
    {\begin{DefinitionBase}}{\begin{DefinitionBase}[#1]}
    \dotfill \, {\em \bf \arabic{claim}} \quad\\{} }
    {\addtocounter{titleclaim}{1}\end{DefinitionBase}}

\newtheorem{TheoremBase}[claim]{Theorem}
\newenvironment{Theorem}[1][-1]{\ifthenelse{\equal{#1}{-1}}
    {\begin{TheoremBase}}{\begin{TheoremBase}[#1]}
    \hrulefill \, {\em \bf \arabic{claim}} \quad\\{} }
    {\addtocounter{titleclaim}{1}\end{TheoremBase}}

\begin{document}

\title{Irreducible Elements in Metric Lattices}
\author{Andreas Lochmann}
\maketitle

\abstract

We describe a natural generalization of irreducibility in order lattices with
arbitrary metrics. We analyse the special cases of valuation metrics and
more general metrics for lattices.

This article is mainly based on a part of the author's doctoral thesis, but
answers some additional questions.

\section{Introduction}

The theory of valuations and metric lattices has been mainly developed and
popularized by John von Neumann and Garrett Birkhoff. In the early years of the
1930s, von Neumann worked on a variation of the ergodic hypothesis, and
inadvertently competed with George David Birkhoff. Only some years later, his
son Garrett Birkhoff pointed von Neumann at the use of lattice theory in Hilbert
spaces. He wrote about this in a note of the Bulletin of the AMS in 1958
\cite{Birkhoff_on_von_Neumann}.

\begin{quote}
John von Neumann's brilliant mind blazed over lattice theory like a meteor,
during a brief period centering around 1935--1937. With the aim of interesting
him in lattices, I had called his attention, in 1933--1934, to the fact that the
sublattice generated by three subspaces of Hilbert space (or any other vector
space) contained 28 subspaces in general, to the analogy between dimension and
measure, and to the characterization of projective geometries as irreducible,
finite-dimensional, complemented modular lattices.

As soon as the relevance of lattices to linear manifolds in Hilbert space was
pointed out, he began to consider how he could use lattices to classify the
factors of operator-algebras. One can get some impression of the initial impact
of lattice concepts on his thinking about this classification problem by reading
the introduction of [...], in which a systematic lattice-theoretic
classification of the different possibilities was initiated. [...]

However, von Neumann was not content with considering lattice theory from the
point of view of such applications alone. With his keen sense for axiomatics, he
quickly also made a series of fundamental contributions to pure lattice theory.
\end{quote}

The modular law in its earliest form (as dimension function) appears in two
papers from 1936 by Glivenko and von Neumann (\cite{Glivenko},
\cite{von_Neumann_continuous_geometries}). Von Neumann used it (and lattice
theory in general) in his paper to define and study Continuous Geometry (aka.
``pointless geometry''), and later applied his knowledge to found Quantum Logic
in his {\it Mathematical Foundations of Quantum Mechanics}. A later survey about
metric posets is \cite{Monjardet}.

The notions of join-irreducibility and join-primeness are fundamental to
Lattice Theory, in the same way as the notion of basis is fundamental to Linear
Algebra (see \cite{Birkhoff}). Hence, it seems plausible to ask for an adaption of
join-irreducibility to metric lattices---the author already used this notion in
\cite{Lochmann_rough_isometries_of_lipschitz_function_spaces} and
\cite{Lochmann_dissertation} to decompose Lipschitz functions and deduce a
rigidity theorem about Lipschitz function spaces. The aim of this article is to
present this new notion of $d$-irreducibility in Section 2 without reference to
Lipschitz function spaces. Section 3 repeats
the definition of a valuation on a lattice and its connection to metrics,
Subsection 3.3 then deduces a characterization of $d$-irreducible elements in
valuation lattices. Subsection 3.2 introduces an alternative definition of
valuation, which is then generalized in Sections 4 and 5 to include further
metrics on lattices, which often are similarly natural but not based on a
valuation. Subsection 5.2 finally deals with the closedness of the subset of
all $d$-irreducible elements in a lattice and in which sense they are a dense
subset of each base.

\subsection{Notation}

Given an element $p$ of a lattice $L$, denote with $\down{p}$ its {\em strictly
lower set}
\beq
\down{p} &:=& \{f\,\in\,L\,:\, f\,<\,p\}.
\eeq
Furthermore, denote with $\Pot (A)$ the power set of $A$.


\section{Irreducibility Relative to a Metric}

Recall the definition of a join-irreducible element $p$
in a lattice $L$:
\beq
p \;=\; f \,\smax\, g &\;\Rightarrow\;& p\; =\; f \quad \n{or} \quad p \;=\; g
\qquad \forall f,\, g\,\in\, L
\eeq
Let $L$ be equipped with the discrete metric $d_\n{dis}$. Then the above
property is equivalent to the following:
\beq
d_\n{dis}\,(p,\, f) \quad\smin\quad d_\n{dis}\,(p,\,g) &\leq& d_\n{dis}\,(p,\,
f\,\smax\,g) \qquad \forall f,\, g\,\in\, L
\eeq
In the same sense, $p$ is completely join-irreducible if and only if
\beq
\biginf_{j\in J}\,d_\n{dis}\,\big(p,\, f_j\big) &\leq& d_\n{dis}\,\left(p,\;
\bigsup_{j\in J} f_j\right) \qquad \forall (f_j)_{j\in J}\subset L,\; J\neq \emptyset.
\eeq

\begin{Definition} \label{DEF_ml-irreducible}
Let $L$ be a lattice with any metric $d$. We call an element $p\,\in\,
L$ {\em $d$-irreducible} if the following holds for all $f,\,g\,\in\,L$:
\beq
d(p,\,f)\;\smin\;d(p,\,g) &\leq& d(p,\, f\smax g)
\eeq
If $L$ is a complete lattice, we call $p$ {\em completely $d$-irreducible}, if
the following holds for all $(f_j)_{j\in J}\subset L$, with $J$ an arbitrary
non-empty index set:
\beq
\biginf_{j\in J}\,d\,\big(p,\, f_j\big) &\leq& d\,\left(p,\;
\bigsup_{j\in J} f_j\right)
\eeq
Denote the subset of $L$ of all completely $d$-irreducible elements with $\mli(L)$.
\end{Definition}

\begin{Proposition} \label{PRO_ml-irred_to_irred}
Let $L$ be a lattice with any metric $d$. Then each $d$-irreducible element is
join-irreducible. However, not every completely $d$-irreducible element
necessarily is completely join-irreducible.
\end{Proposition}
\begin{Proof}
Let $p\,\in\,L$ be $d$-irreducible and $p\,=\,f\smax g$. Then $d(p,\, f\smax
g)\,=\,0$ and hence either $d(p,\,f)\,=\,0$ or $d(p,\,g)\,=\,0$ (or both).

For a counter-example
to complete join-irreducibility, let $L\,=\,[0,\,1]$ with standard metric,
supremum and infimum. Take $f_n\,=\, 1- 1/n$, $n\,\in\, \N^*$, then
$p\,=\,1\,=\,\bigsup f_n$, hence $p$ is not completely join-irreducible. Still,
it is completely $d$-irreducible: Any sequence of real numbers $f_n$ with
$p\,=\, \bigsup f_n$ must converge to $p$ from below, hence $\biginf d(p,\,
f_n)\, =\, 0$.
\end{Proof}

As a consequence, if $L$ is a complemented lattice, join-irreducibility,
complete join-irreducibility, $d$-irreducibility, and complete
$d$-irreducibility are all equivalent; the irreducible elements are simply
those with trivial strictly lower set.


\section{Valuations}

\begin{Definition} \label{DEF_valuation}
A {\em valuation} on a lattice $L$ is a function $v: L\rightarrow \R$ which
satisfies the {\em modular law}
\beq
v(f) \;+\; v(g) &=& v(f\smin g) \;+\; v(f\smax g) \quad \forall f,g \in L.
\eeq
A valuation $v$ on $L$ is called {\em isotone} [{\em positive}] if for all
$f,g\in L$ the relation $f < g$ implies $v(f) \leq v(g)$ [$v(f) < v(g)$].
\end{Definition}

If $L$ is totally ordered, then each function $v: L \rightarrow \R$ is a
valuation. It is isotone [positive] if and only if $v$ is [strictly]
monotonically increasing.

Valuations can be used to define metrics on lattices, as the following Lemma
demonstrates. It is a part of Theorem X.1 and a note in subsection X.2 of
\cite{Birkhoff}, and is proved there. An alternative proof is given in
\cite{Lochmann_dissertation}.

\begin{Lemma} \label{LEM_birkhoff_metric_lattice}
Let $v$ be an isotone valuation on the distributive lattice $L$. Then 
\beq
d_v(f,g) &:=& v(f\smax g) \;-\; v(f\smin g)
\eeq
defines a pseudo-metric with the following properties:
\begin{enumerate}
\item
If there is a least element $0\in L$, then
\beq
v(f) &=& v(0) \;+\; d_v(f,0) \quad \n{for all} \quad f \in L,
\eeq
\item
$d_v$ is a metric if and only if $v$ is positive.
\end{enumerate}
We call $d_v$ a {\em valuation (pseudo-)metric}. A lattice together with a
valuation metric is sometimes called a {\em metric lattice}; however, as we will
deal with lattices with non-valuation metrics as well (particularly the
supremum metric), we should better distinguish between {\em valuation metric
lattices} and {\em non-valuation metric lattices}.
\end{Lemma}

\subsection{Examples}

Valuations and valuation metrics arise in a multitude of situations:

\begin{Example}
Let $L=(\N^*, \gcd, \lcm)$. Then each logarithm is a positive valuation on
$L$. The join-irreducible, completely join-irreducible, $d$-irreducible and
completely $d$-irreducible elements are exactly the prime powers.
\end{Example}

\begin{Example}
Let $(X, \Sigma, \my)$ be a probability space. The $\sigma$-algebra $\Sigma$ is
a Boolean lattice by union and intersection. Let $c\in \R$ be arbitrary, then
\beq
v(A) &:=& \my(A) + c
\eeq
defines an isotone valuation on $\Sigma$ with $v(\emptyset) = c$. The valuation
$v$ is positive if and only if there are no null sets in $X$ other than
$\emptyset$.
%
%
The distance function $d_v(A,B) := v(A\cup B) - v(A\cap B)$ is the measure of the
symmetric difference $A\sd B$ of $A$ and $B$, if $A\sd B\in \Sigma$. It relates
to the Hausdorff distance just as the 1-distance of functions relates to the
supremum distance.
\end{Example}

\begin{Example} \label{EXA_PGV_metric}
Let $V$ be any finite dimensional vector space, and $L=\PG(V)$ its lattice of
subvector spaces, with $\smin$ the intersection and $\smax$ the span (the {\em
projective geometry} of $V$). Then the dimension function is a positive
valuation on $L$. (This is the similarity between dimension and measure
mentioned before in \cite{Birkhoff_on_von_Neumann}.) The join-irreducible,
completely join-irreducible, $d$-irreducible and completely $d$-irreducible
elements are exactly the one-dimensional subspaces and the zero-dimensional one.
\end{Example}

\begin{Example} \label{EXA_l1_metric}
Let $X$ be a measure space and $L$ the space of integrable Lipschitz functions of
Lipschitz constant $\leq\,1$. We may apply the Lebesgue integral to
gain an isotone valuation on $L$; as $f + g = (f\smin g) + (f\smax g)$ holds
pointwise, we conclude
\beq
\int f\, \d \my + \int g\, \d \my &=& \int (f\smin g)\, \d \my + \int (f\smax
g)\, \d \my.
\eeq
If $X$ is a Euclidean space, or a discrete space without non-trivial null sets,
this valuation is positive, because any non-trivial non-negative Lipschitz
function has positive Lebesgue integral. Positivity fails in cases
where $X$ contains an isolated point or continuum of measure zero.

As $|f-g| = (f\smax g) - (f\smin g)$ holds pointwise, the valuation metric $d_v$
equals the $L^1$-distance defined by
\beq
d_1(f,g) &:=& \int |f-g|\,\d\my.
\eeq

Each function $\Lambda(x,\,r):\,L\,\rightarrow\,[0,\infty)$ of the form
\beq
\Lambda(x,\,r)(y) &:=& 0\,\smax\,\big(r\,-\,d(x,\,y)\big)
\eeq
with $x\,\in\,X$ and $r\,\in\,[0,\infty)$ is join-irreducible, but not
necessarily completely join-irreducible. In general, the only $d$-irreducible
function is the zero function.

The $L^1$-metric can be slightly modified to yield other valuation metrics: Let
$\kappa: [0,\infty) \rightarrow [0,\infty)$ be a positive valuation (i.e.,
strictly monotonically increasing), then
\beq
v_{\my,\kappa}(f) &:=& \int \kappa(f(x))\,\d \my(x)
\eeq
is a positive valuation.
\end{Example}


\subsection{Difference Valuations} \label{SEC_difference_valuations}

A nearly equivalent approach to valuations is to use difference valuations:

\begin{Definition} \label{DEF_difference_valuation}
A {\em difference valuation} on a distributive lattice $L$ is a function $w:
L\times L \rightarrow \R$ which satisfies the cut law
\beq
w(f,\,g) &=& w(f,\, g\smax h) \;+\; w(f\smin h,\, g).
\eeq
A difference valuation $w$ is called {\em isotone} if its values are
non-negative, and {\em positive}, if $w(f,\,g)\,=\,0$ implies $f\,\leq\,g$.
\end{Definition}

\begin{figure}[t]
\centerline{\epsfig{file=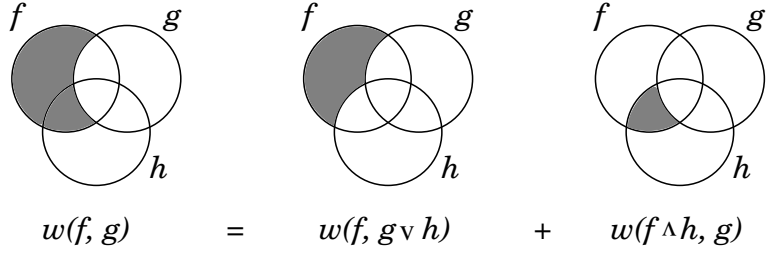, width = 10cm}}
\caption{Visualization of the cut law of difference valuations using Venn diagrams. Given
a Stone representation $\pi$, the set $\pi(f)\setminus \pi(g)$ is cut along $\pi(h)$ to
give $\pi(f)\setminus \pi(g\,\smax\, h)$ and $\pi(f\,\smin\, h)\setminus \pi(g)$.}
\label{FIG_cut_law}
\end{figure}

Given a valuation $v$ on a distributive lattice $L$,
\beq
w(f,\,g) &:=& v(f) - v(f\smin g)
\eeq
defines a difference valuation, as one can easily check. The cut law follows
from the modular equality and vice versa---it has been dubbed ``cut law''
because of its appearance when applied to sets in a Venn diagram, see Figure
\ref{FIG_cut_law}. The difference valuation is isotone/positive
if and only if the valuation $v$ is isotone/positive. If $L$ admits a least
element $0$, each [isotone/positive] difference valuation $w$ in turn defines an
[isotone/positive] valuation $v$ by 
\beq
v(f) &:=& w(f,\,0) \;+\; c
\eeq
for any $c\,\in\,\R$, and any valuation of $L$ with difference valuation $w$ is
of this form.
Finally, the distance function $d$ of a valuation can be equally well expressed
as
\beq
d(f,\,g) &=& w(f,\,g) \;+\; w(g,\,f).
\eeq

\begin{figure}[t]
\centerline{\epsfig{file=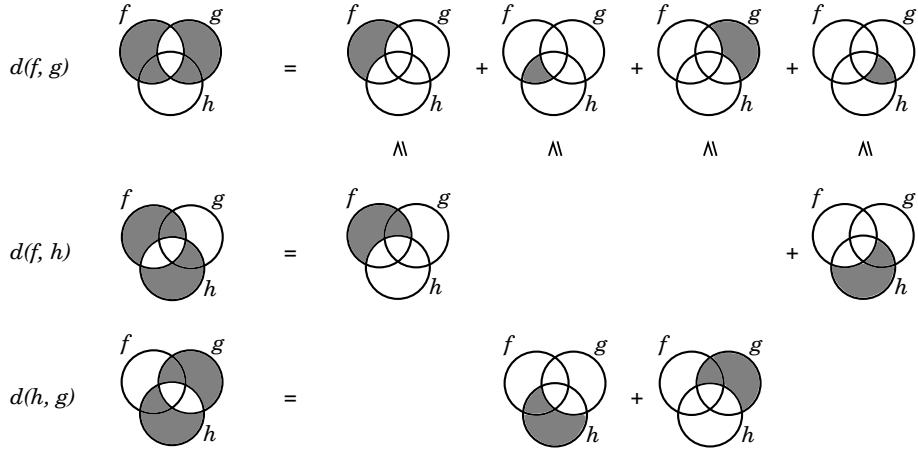, width = 12cm}}
\caption{Proof of the triangle inequality for valuation metric lattices using
difference valuations and Venn diagrams. Note that $f, g, h$ are elements of an
arbitrary distributive lattice, and represented by sets via Stone duality.}
\label{FIG_dv_triangle_inequality}
\end{figure}

If the lattice $L$ is complemented, $w(f,\,g)$ equals $v(f\setminus g)$.
Difference valuations are easier to use in cases where a lattice is not
complemented, as they can be used as substitutes for the relative complement
operation in calculations with metrics. For example, the proof of Lemma
\ref{LEM_birkhoff_metric_lattice} can be seen by a simple application of Venn
diagrams (see Figure \ref{FIG_dv_triangle_inequality}); for details and further
examples to deduce metric inequalities in order lattices see \cite{Lochmann_dissertation}.

\subsection{$d$-Irreducible Elements}

As a triviality, in the definition of a join-irreducible element,
\beq
p \;=\; f \,\smax\, g &\;\Rightarrow\;& p\; =\; f \quad \n{or} \quad p \;=\; g
\qquad \forall f,\, g\,\in\, L,
\eeq
the elements $f$ and $g$ may be chosen to be $\in\,\down{p}\,\subset\,L$. This
accounts for $d$-irreducible elements as well, but is less trivial:

\begin{Lemma} \label{LEM_irreducible_just_down_valuation}
Let $L$ be a distributive lattice, and $d$ a positive valuation metric on $L$.
$p\,\in\,L$ is $d$-irreducible if and only if
\beq
d(p,\,f)\;\smin\;d(p,\,g) &\leq& d(p,\, f\smax g)
\eeq
holds for all $f,\,g\,\in\down{p}$. In this case, ``$\leq$'' can be replaced by ``$=$''.

If $L$ is completely distributive, then
the analog holds for complete $d$-irreducibility as well.
\end{Lemma}
\begin{Proof}
Let $f,\,g\,\in\,L$ be arbitrary and $p\,\in\,L$ as above. Then holds:
\beq
d(f\,\smax\,g,\,p) &=& w(p,\,f\,\smax\, g) \;+\; w(f\,\smax\,g,\,p) \\
&\geq& w(p,\,f\,\smax\,g) \;+\; \big(w(f,\,p) \,+\,w(g,\,p)\big) \\
&=& d\big((f\smax g)\,\smin\, p,\, p\big) \;+\; \big(w(f,\,p) \,+\,w(g,\,p)\big) \\
&\geq& \big(d(f\,\smax\,p,\,p)\,\smin\, d(g\,\smax\,p,\,p)\big) 
       \;+\; \big(w(f,\,p) \,+\,w(g,\,p)\big) \\
&=& \big(w(p,\,f)\,\smin\,w(p,\,g)\big) \;+\; w(f,\,p) \,+\,w(g,\,p) \\
&\geq& \big(w(p,\,f)\,+\,w(f,\,p)\big) \;\smin\; \big(w(p,\,g)\,+\,w(g,\,p)\big) \\
&=& d(f,\,p) \;\smin\; d(g,\,p)
\eeq
(1: definition, 2: by cut law, 3: definition, 4: hypothesis, 5: definition, 6:
by distributivity and positivity of $w$, 7: definition). Each step holds in the infinite case as
well, one only has to add in step 6 that $L$ is completely distributive.

For equality, note that $d(p,\,f) \,\geq\,d(p,\,f\smax g)$ is obvious because
$f\,\leq\,f\smax g\,\leq p$; same holds for $g$ and thus 
\beq
d(p,\,f)\;\smin\;d(p,\,g) &\geq& d(p,\, f\smax g).
\eeq
\end{Proof}

There is a characterization of join-irreducibility of an element $p\,\in\,L$ in
terms of its strictly lower set $\down{p}$: $p$ is join-irreducible if and only if
for each $f,\,g\,\in\down{p}$ holds $f\,\smax\,g\,\in\down{p}$, i.e. if and
only if $\down{p}$ is join-closed. Analo\-gously, $p$ is a completely
join-irreducible element of a complete lattice $L$ if and only if $\down{p}$ is
join-complete (i.e. each supremum 
of elements of $\down{p}$ again is contained in $\down{p}$).
For valuation metrics, there is a similar characterization of
$d$-irreducibility:

\begin{Theorem} \label{THE_d-irreducibility_for_valuation_metrics}
Let $L$ be a distributive lattice, and $d$ a positive valuation metric on $L$.
An element $p\,\in\,L$ is $d$-irreducible if and only if the strictly
lower set $\down{p}$ is totally ordered.
%
%
\end{Theorem}
\begin{Proof}
``$\Rightarrow$'': Let $f,\,g\,\in\down{p}$ be arbitrary.
\beq
d(f\,\smax\,g,\,p) &=& d(f,\,p) \;\smin\; d(g,\,p)\\
&=& w(p,\,f) \;\smin\; w(p,\,g) \\
&=& \big(w(g,\,f) \,+\, w(p,\, f\smax g)\big)
    \;\smin\; \big(w(f,\,g) \,+\, w(p,\, g\smax f)\big) \\
&=& \big(w(g,\,f) \,\smin\, w(f,\,g)\big) \;+\; w(p,\, g\smax f) \\
&=& \big(w(g,\,f) \,\smin\, w(f,\,g)\big) \;+\; d(f\,\smax\, g,\, p)
\eeq
and hence $w(g,\,f)\,\smin\, w(f,\,g)\,=\,0$, thus one of them is zero, and we
have either $f\,\leq\,g$ or $g\,\leq\,f$.

``$\Leftarrow$'': Let $f,\,g\,\in\down{p}$ be arbitrary (see Lemma
\ref{LEM_irreducible_just_down_valuation} why we may restrict to $\down p$). As
$\down{p}$ is totally ordered, $f\smax g$ is $f$ or $g$, and hence the
condition for $d$-irreducibility is trivial.

%

\end{Proof}

Theorem \ref{THE_d-irreducibility_for_valuation_metrics} shows that
$d$-irreducibility does not depend on the concrete choice of a valuation metric
for the lattice $L$. This result resembles an earlier connection found in
Lipschitz function spaces: If $L$ is the space of bounded non-negative Lipschitz functions
of a metric space $X$ with Lipschitz constant $\leq 1$ with pointwise supremum
and infimum and supremum metric $d_\infty$, then the completely
$d_\infty$-irreducible elements are exactly those functions of the form
\beq
\Lambda(x,\,r): \;\; L &\rightarrow& [0,\,\infty)\\
y &\mapsto& \big(r\,-\,d_X(x,\,y)\big)\;\smax\;0
\eeq
with $x\,\in\,L$ and $r\,\in\,[0,\,\infty)$ (see Example \ref{EXA_l1_metric},
\cite{Lochmann_rough_isometries_of_lipschitz_function_spaces}, \cite{Lochmann_dissertation}).
Although the supremum metric $d_\infty$ is not a valuation metric, but an
intervaluation metric (see Definition \ref{DEF_intervaluation}), its completely
$d_\infty$-irreducible elements are fully determined without any reference to
the chosen metric on $L$. One might even get rid of the metric on $X$ by
referring to minimal functions with a given function value at a single point.


\section{Ultravaluations}

One advantage of the definition of difference valuations in Subsection
\ref{SEC_difference_valuations} is the following alternative to valuations in
lattices, which adds further examples to our list of metrics on lattices and is
easily described in terms of a variant of Definition \ref{DEF_difference_valuation}.

\begin{Lemma} \label{LEM_ultravaluation_metric}
Let $L$ be a distributive lattice, and let $w:L\times L \rightarrow [0,\infty)$
be a map which satisfies
\beq
\n{(1)} && w(f,g) \;\;=\;\; 0 \quad\n{whenever}\quad f\leq g,\quad\n{and}\\
\n{(2)} && w(f,\,g) \;\;=\;\; w(f\smin h,\, g) \;\smax\; w(f,\, g\smax h) \quad
\forall f,g,h\in L.
\eeq
We call $w$ a {\em difference ultravaluation}, or just {\em ultravaluation}.
Define
\beq
d_w(f,g) &:=& w(f,g) \;\smax\; w(g,f).
\eeq
Then $d_w$ is a
pseudo-ultrametric. $d_w$ is an ultra\imetric if and only if
$w(f,g) = 0 \,\Rightarrow\, f\leq g$ holds.
\end{Lemma}
\begin{Proof}
To get from difference valuations to ultravaluations, we just replaced all
occurences of ``$+$'' by ``$\smax$''. As both operations are associative and
commutative, we can transfer most proofs of valuations just by replacing
``$+$'' by ``$\smax$'', this includes the proof of the triangle inequality:
\beq
d_w(f,g) &=& w(f, g\smax h) \smax w(f\smin h, g) \smax w(g, f\smax h) \smax
w(g\smin h, f)\\
w(f, g\smax h) &\leq& w(f, h) \quad \n{etc.}\\
\Rightarrow \quad d_w(f, g) &\leq& w(f, h) \smax w(h,g) \smax w(g,h) \smax w(h,
f) \;=\; d_w(f, h) \smax d_w(h, g)
\eeq

On the other hand, contrary to the valuation case, the property $d_v(f,f) = 0$
does not follow from property (2) -- we have to conclude it from (1).

Assume $w(f,g)=0 \,\Rightarrow\, f\leq g$ holds. Let $d_w(f,g) = 0$. This
implies $w(f,g) = 0$ and $w(g,f) = 0$, and hence $f\leq g$, $g\leq f$, and
$f=g$. Now assume $d_w$ is a metric, $f\nleq g$, and $w(f,g) = 0$. Then
\beq
w(f, f\smin g) \;=\; w(f\smin g, f\smin g) \;\smax\; w(f,g) \;=\; 0.
\eeq
Due to $f\nleq g$, we have $f\neq f\smin g$, hence
\beq
0 \;<\; d_w(f, f\smin g) \;=\; w(f,f\smin g) \;\smax\; w(f\smin g, f) \;=\;
w(f\smin g, f).
\eeq
But $f\smin g \leq f$, contradiction. 
\end{Proof}

\subsection{Examples}

\begin{Example} \label{EXA_ultravaluation_by_function}
Let $X$ be any set, $\kappa:X\rightarrow [0,\infty)$ arbitrary and fixed, and $L$ a
lattice of subsets of $X$. For $A,B\in L$ consider
\beq
w(A,B) &:=& 0\;\smax\; \sup_{x\in A\setminus B}\; \kappa(x).
\eeq
$w$ defines an ultravaluation.

Choose $\kappa$ to be a positive constant, then the ultrametric resulting from
$w$ will be the discrete metric on $L$.
\end{Example}

\begin{Example} \label{EXA_ultravaluation_lip}
Let $X$ be any metric space and $\Lip_0 X$ its lattice of bounded Lipschitz
function of Lipschitz constant $\leq 1$. Besides its Stone representation, we
want to provide another, more intuitive
representation of the space $\Lip_0 X$ by a lattice of sets, using its {\em hypograph}
(cp. ``epigraph'' in \cite{Rockafellar})
\beq
\hyp: \;\Lip X &\rightarrow& \Pot\,(X\times [0,\infty)) \\
f &\mapsto& \{(x,r)\;:\; f(x) \leq r\}.
\eeq
$(\im \hyp,\, \cap,\, \cup)$ obviously is isomorphic to $(\Lip_0 X,\,\smin,\smax)$
as a lattice; however, they are not yet isomorphic as {\em complete} lattices:
Infinite unions of the closed sets in $\im\hyp$ are not closed in general -- we
have to use the union with closure ``$\ \bar{\cup}$'' instead of the traditional
union.  (Alternatively, we could identify subsets of $X\times [0,\infty)$ with
the same closure.)

We now apply Example \ref{EXA_ultravaluation_by_function}. The most canonical
$\kappa$ would be $\kappa = \pi_2$, the projection onto $[0,\infty)$. The
corresponding ultrametric on $L$ is
\beq
d_\kappa(f,g) &=& 
0\;\smax\; \sup\;\{f(x) \,\smax\, g(x)\n{ with $x\in X$ such that } f(x) \neq
g(x)\}.
\eeq
We shall call this metric the ``peak metric'' on $\Lip X$.

Another possible choice for $\kappa$ is as follows: Choose a basepoint $x_0\in
X$ and $\kappa(x,r) \;:=\; d_X(x, x_0)$. Then $d_\kappa$ will describe the
greatest distance from $x_0$ at which $f$ and $g$ still differ. Finally,
$\kappa(x,r) \;:=\; \exp(-d_X(x,x_0))$ will describe the least
distance from $x_0$ at which $f$ and $g$ differ. We will call the first case the
``outer basepoint metric'' and the second case the ``inner basepoint metric''.

An application of the lower basepoint metric is as follows: Given a free group
$F$ with neutral element $x_0$, identify each normal subgroup $N\trianglelefteq
F$ with its characteristic function on $F$. These are 1-Lipschitz functions in
the canonical word metric of $F$. $d_\kappa$ then defines a topology on $\Lip
F$, which restricts to the Cayley topology (\cite{dlHarpe}, V.10) on the subset
of normal subgroups.

The $\Lambda$-functions defined in Example \ref{EXA_l1_metric} are exactly the
$d$-irreducible functions of the peak metric. The $d$-irreducible functions of
the outer basepoint metric are those functions $\Lambda(x,\,r)$ with
$x\,\neq\,x_0$, the inner basepoint metric doesn't admit any non-trivial
$d$-irreducible function in general. Finally, none of these three metrics admits a
non-trivial completely $d$-irreducible function.
\end{Example}

\begin{Lemma} \label{LEM_finite_ultravaluation_lattice}
Let $X$ be finite, and let $L$ be a lattice of subsets of $X$. Then any ultravaluation on
$L$ is of the form of Example \ref{EXA_ultravaluation_by_function}.
\end{Lemma}
\begin{Proof}
For $x\,\in\, X$ and $A,\,B\,\subset\,X$ define
\beq
\kappa(x) &:=& \inf\;\{w(C,\,D)\;:\; C,\,D\,\in\, L\;\n{with}\; x\ \in\  C,\; x\notin D\}\\
\n{and} \quad w'(A,\,B) &:=& 0\;\smax\; \sup_{y\in A\setminus B}\; \kappa(y).
\eeq
Assume $w'(A,\,B)\,>\,w(A,\,B)$. Then there is $y\,\in\,A\setminus B$ with
$\kappa(y) \,\geq\, w(A,\,B)$, but this cannot happen, as one may choose $C=A$
and $D=B$. Hence, assume $w'(A,\, B)\,<\,w(A,\,B)$. Then for all $y\,\in\,
A\setminus B$ there should be $C,\,D\,\in\,L$ with $y\,\in\,C\setminus D$ and
$w(C,\,D) \,<\, w(A,\, B)$. As
\beq
w(C,\,D) \,\geq\, w(C\,\smin A,\, D\,\smax\, B),
\eeq
we might choose without loss of generality $C\subset A$ and $D\supset B$, as
choosing $C\cap A$ instead of $C$ and $D\cup B$ instead of $D$ further decreases
$w(C,\,D)$. The cut law now yields
\beq
w(A,\,B) &=& w(C\smin D,\, B) \,\smax\, w(C,\, D)\,\smax\, w(A\smax D,\, B\smax
C)\,\smax\, w(A,\, C\smax D).
\eeq
As $w(C,\, D)\,<\, w(A,\,B)$ by assumption, we find that at least one of $(C\cap
D)\setminus B$, $(A\cup D)\setminus (B\cup C)$, and $A\setminus(C\cup D)$ must
be non-empty. Choose $y'$ out of their union and repeat the above argument for
the now smaller subset. We get an infinite sequence of different elements from
$X$, which is a contradiction because $X$ is finite.
\end{Proof}

\begin{Example} \label{EXA_hausdorffdimension}
Not all ultravaluations are of the kind of Example
\ref{EXA_ultravaluation_by_function}. Let $X$ be any metric space and $L$ the
lattice of subsets of $X$. Define $w(A,\,B)$ to be the Hausdorff dimension
of $A\setminus B\,\in\,L$ plus $1$, and $0$ if $A\setminus B\,=\,\emptyset$.
Then $w$ is an ultravaluation and $d_w$ an ultrametric.

The $d$-irreducible subsets and the completely $d$-irreducible subsets are
exactly the join-irreducible subsets, namely those with one or zero elements,
because $L$ is complemented.
\end{Example}

Comparing Examples \ref{EXA_PGV_metric} and \ref{EXA_hausdorffdimension}, one
should note that the join operation in the former is the span, but in the latter
is the union. Thus, the first example gives rise to a valuation, the second one
to an ultravaluation.


\subsection{$d$-Irreducible Elements} \label{SEC_ultravaluations_irreducibles}

Lemma \ref{LEM_irreducible_just_down_valuation} can be easily adapted to the
case of ultravaluations by replacing all remaining ``$+$'' by ``$\smax$''.
Indeed, Lemma \ref{LEM_irreducible_just_down_valuation} holds in an even
broader generalization, what we will demonstrate in Lemma
\ref{LEM_irreducible_just_down_intervaluation}.

When following the proof of Theorem
\ref{THE_d-irreducibility_for_valuation_metrics} for ultravaluation metrics
(remember that join-irreducibility is $d_\n{dis}$-irreducibility for the
discrete metric $d_\n{dis}$, which is an ultravaluation metric), one ends up
with the following inequality:
\beq
d(f,\,g)&\leq& d(p,\,f\,\smax\,g)
\eeq
for all $d$-irreducible elements $p$ and all $f,\,g\,\in\down{p}$. If $L$
contains a least element $0\,\in\,L$, we conclude as special case
\beq
d(0,\,g) &\leq& d(g,\,p) \qquad \forall g\,\in\down{p}.
\eeq

One would hope that there is a similar characterization of $d$-irreducible
elements in the ultravaluation case as it is in the valuation case. Starting
from the case of the discrete metric, one would ask whether join-irreducibility
is exactly this characterization, i.e. whether all join-irreducible elements
are $d$-irreducible for any ultravaluation metric $d$. This, however, is wrong.

\begin{Example}
We refer to Example \ref{EXA_ultravaluation_by_function}.
Let $X \,=\, \{1,\,2,\,3\}\,\subset\,\Z$ and let $\kappa$ be the identity. Let
$L$ be the lattice $\{\emptyset,\,\{2\},\,\{3\},\{2,\,3\},\,X\}$ of subsets of $X$. Then
$X\,\in\,L$ is join-irreducible (because it is the only set containing $1$), but not
$d$-irreducible: $d(X,\, \{2\})\,=\, 3$, $d(X,\,\{3\})\,=\,2$ and $d(X,\,
\{2,\,3\})\,=\,1$. In particular, this example shows that $d$-irreducibility
depends on the concrete choice of $\kappa$, respectively on the choice of the
ultravaluation.
\end{Example}

\begin{Question}
Is there a nice criterion to decide whether all join-irreducible elements in an
ultravaluation metric lattice are $d$-irreducible?
\end{Question}

\smallskip

\begin{figure}[t]
\centerline{\epsfig{file=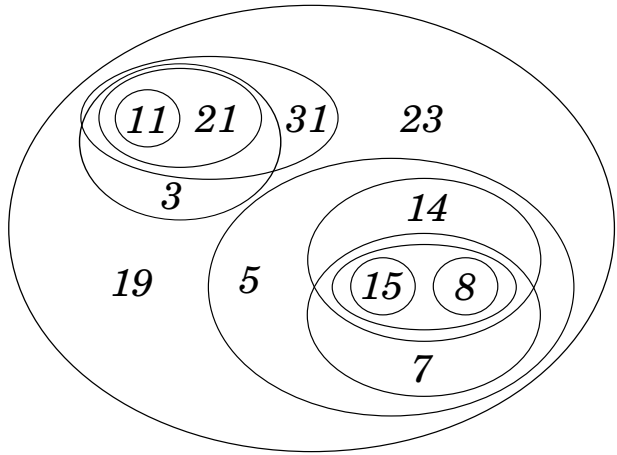, width = 8cm}}
\caption{We refer to Example \ref{EXA_ultravaluation_by_function}. Let $L$ be
the lattice of sets spanned by the shown sets of natural numbers and let
$\kappa$ be the identity. A set $A$ is {\em not} $d$-irreducible, if and only if
there are subsets $B$ and $C\,\in\,L$ of $A$, such that both $B$ and $C$
contain at least one number each, which is larger than any of the remaining
numbers in $A\setminus(B\,\cup\,C)$.
Which of the shown subsets are $d$-irreducible?} 
\label{FIG_irreducibility_puzzle}
\end{figure}

Lemma \ref{LEM_finite_ultravaluation_lattice} characterizes all finite
ultravaluation lattices. However, finding the $d$-irreducible subsets in a
finite ultravaluation lattice can still be non-trivial. We demonstrate this by
restating the problem as a puzzle in Figure \ref{FIG_irreducibility_puzzle}
and leave it to the reader to find any patterns.

\section{Intervaluations and Topological Aspects}

We now present a generalized notion of valuation which includes normal
valuations and ultravaluations. In addition, this notion of
intervaluations also includes the supremum metric of function spaces, just as
the $L^1$-metric was found to be a valuation in Example \ref{EXA_l1_metric}.

Similar to the case of the ultravaluation, we first recognize the possibility to
replace ``$+$'' in the definition of a difference valuation by any commutative
and associative binary operation. But this alone will not suffice to encompass
the supremum metric, we have to weaken the main property of a difference
evaluation as well:

\begin{Definition} \label{DEF_intervaluation}
An {\em intervaluation} on a distributive lattice $(L,\smin,\smax)$ is a
map $w: L\rightarrow [0,\infty)$ together with a commutative and associative
binary operation $\circ_w: [0,\infty) \times [0,\infty) \rightarrow [0,\infty)$,
such that the following properties hold:
\begin{enumerate}
\item $r\,\circ_w\, 0 \;=\; 0\,\circ_w\, r \;=\; r$
\item $r\,\circ_w\, t \;\leq\; (r \,+\, s)\,\circ_w\, (t \,+\, u) \;\leq\;
(r\,\circ_w\,t) \,+\, (s\,\circ_w\,u)$
\item $r\,\smax\, s \;\leq\; r\, \circ_w\, s$ \quad\n{(follows from (1) and (2))}
\item $f\, \leq\, g \quad\Rightarrow\quad w(f,\,g)\, =\, 0$
\item $w(f,\, g\smax h)\, \circ_w\, w(f\smin h,\, g) \;\leq\; w(f,\, g)
\;\leq\; w(f,\, g\smax h) \,+\, w(f\smin h,\, g)$ \\ (left and right modular
inequality, or cut law)
\end{enumerate}
for all $f,g,h\in L$ and $r,s,t,u\in [0,\infty)$. The corresponding {\em
intervaluation metric} then is defined to be
\beq
d_w(f,\,g) &:=& w(f,\,g)\,\circ_w\,w(g,\,f).
\eeq
The intervaluation is {\em positive} if
\beq
w(f,\,g) \,=\,0 &\Rightarrow& f \,\leq\,g.
\eeq
\end{Definition}

\begin{Proposition} \label{PRO_intervaluation}
An intervaluation $w$ on $L$ and its metric $d_w$ always fulfill:
\begin{enumerate}
\item $w(f,\,g) \;=\; w(f\smax g,\,g) \;=\; w(f,\,f\smin g)\;=\; d_w(f\smax
g,\, g)\quad \forall f,\,g\in L$.
\item $d_w$ is a pseudo-\imetric.
\item $d_w$ is a \imetric if and only if $w$ is positive.
\end{enumerate}
\end{Proposition}
\begin{Proof}
{\bf (1)} We choose $h=f$ or $h=g$ in both modular inequalities:
\beq
&&0\;\circ_w\; w(f,\,g) \quad\leq\quad w(f\smax g,\,g) \quad\leq\quad 0 \;+\;
w(f,\,g)\\
&& w(f,\,g) \;\circ_w\; 0 \quad\leq\quad w(f,\,g\smin f) \quad\leq\quad w(f,\,g)
\;+\; 0\\
\n{and} &&d_w(f\smax g,\, g) \quad=\quad w(f\smax g,\, g)\;\circ_w\; 0
\quad=\quad w(f,\, g).
\eeq

{\bf (2)} From the definition we see $d_w(f,\, g)\geq 0$ and $d_w(f,\, f) =
0$ for all $f,\,g\in L$. As $\circ_w$ is commutative, $d_w$ is symmetric.
\beq
d_w(f,g) &=& w(f,\,g) \;\circ_w\; w(g,\,f)\\
&\leq& \left(w\left(f\smin h,\, g\right) \;+\;
w\left(f,\,g\smax h\right)\right) \;\circ_w\; \left(w\left(g\smin h,\,
f\right) \;+\; w\left(g,\,f\smax h\right)\right)\\
&\leq& \left(w\left(h,\,g\right) \;+\; w\left(f,\,h\right)\right)\;\circ_w\;
\left(w\left(h,\,f\right) \;+\; w\left(g,\,h\right)\right))\\
&=& \left(w\left(f,\,h\right) \;+\; w\left(h,\,g\right)\right)\;\circ_w\;
\left(w\left(h,\,f\right) \;+\; w\left(g,\,h\right)\right))\\
&\leq& \left(w\left(f,\,h\right) \;\circ_w\; w\left(h,\,f\right)\right)\;+\;
\left(w\left(h,\,g\right) \;\circ_w\; w\left(g,\,h\right)\right))\\
&=& d_w(f,\, h) \;+\; d_w(h,\, g)
\eeq

{\bf (3, ``$\Rightarrow$'')} Assume $0 = w(f,\,g) = w(f,\,f\smin g)$. Then
$d_w(f,f\smin g) = 0 + 0 = 0$. As $d_w$ is a metric, we have $f=f\smin g$, so
$f\leq g$.

{\bf (3, ``$\Leftarrow$'')} $d_w(f,\,g) = 0$ implies $w(f,\,g) = 0$ and
$w(g,\,f)=0$, hence $f\leq g\leq f$, and $f=g$.
\end{Proof}

We now show the generalization of Lemma
\ref{LEM_irreducible_just_down_valuation} for intervaluations, which we already
announced in subsection \ref{SEC_ultravaluations_irreducibles}.

\begin{Lemma} \label{LEM_irreducible_just_down_intervaluation}
Let $L$ be a distributive lattice, and $d$ a positive intervaluation metric on $L$.
$p\,\in\,L$ is $d$-irreducible if and only if
\beq
d(p,\,f)\;\smin\;d(p,\,g) &\leq& d(p,\, f\smax g)
\eeq
holds for all $f,\,g\,\in\down{p}$. In this case, ``$\leq$'' can be replaced by ``$=$''.

If $L$ is completely distributive, then
the analog holds for complete $d$-irreducibility as well.
\end{Lemma}
\begin{Proof}
Let $f,\,g\,\in\,L$ be arbitrary and $p\,\in\,L$ as above. Then holds:
\beq
d(f\,\smax\,g,\,p) &=& w(p,\,f\,\smax\, g) \;\circ_w\; w(f\,\smax\,g,\,p) \\
&\geq& w(p,\,f\,\smax\,g) \;\circ_w\; \big(w(f,\,p) \,\circ_w\,w(g,\,p)\big) \\
&=& d\big((f\smax g)\,\smin\, p,\, p\big) \;\circ_w\; w(f,\,p) \,\circ_w\,w(g,\,p) \\
&\geq& \big(d(f\,\smax\,p,\,p)\,\smin\, d(g\,\smax\,p,\,p)\big)
       \;\circ_w\; w(f,\,p) \,\circ_w\,w(g,\,p) \\
&=& \big(w(p,\,f)\,\smin\,w(p,\,g)\big) \;\circ_w\; w(f,\,p) \,\circ_w\,w(g,\,p) \\
&\geq& \big(w(p,\,f)\,\circ_w\,w(f,\,p)\big) \;\smin\; \big(w(p,\,g)\,\circ_w\,w(g,\,p)\big) \\
&=& d(f,\,p) \;\smin\; d(g,\,p)
\eeq
(1: definition, 2: by left modular inequality, 3: definition, 4: hypothesis, 5: definition, 6:
by cases and monotony of ``$\circ_w$'' (property (2) in Definition \ref{DEF_intervaluation}), 7:
definition). Each step holds in the infinite case as well.
\end{Proof}

\subsection{Examples}

\begin{Example} \label{EXA_intervaluation_operations}
There are several possible choices for the commutative and associative binary
operation $\circ_w$ in Definition \ref{DEF_intervaluation}. Choosing addition
leads directly to the definition of valuations. The next important choice is the
maximum operation: Properties (1) and (3) are obviously fulfilled, the left
side of (2) as well. (2.right) needs some short consideration: As
$+$ distributes over $\smax$, the right-hand side equals
\beq
(r\,\smax\,t)\,+\,(s\,\smax\,u) &=&
(r+s)\,\smax\,(r+u)\,\smax\,(t+s)\,\smax\,(t+u)
\eeq
which is greater or equal $(r+s)\,\smax\,(t+u)$ for all $r,\,s,\,t,\,u\,\in
[0,\,\infty)$.

Each norm $||\cdot||$ on $\R^2$ with certain normalization properties qualifies
as an operation $\circ_w$ via $r\,\circ_w\,s\,:=\, ||(r,s)||$. This accounts for
the $\ell^p$-norms:
\beq
r\,\circ_p\,s \;:=\; \big|\big|(r,\,s)\big|\big|_p\;:=\;\sqrt[p]{r^p\,+\,s^p}
\eeq
for $p\,\in\, [1,\infty)$. Again, properties (1), (2.left) and (3) of
Definition \ref{DEF_intervaluation} are trivial. Property (2.right) is the
triangle inequality of the $\ell^p$-norms (i.e.\ a special case of the Minkowski
inequality \cite{Werner}).
\end{Example}

Given any metric $d$ on $L$ we may define $w_d(f,\,g) \;:=\; d(f\smax g,\, g)$
and deduce $\circ_w$ from $d(f,\,g) = w_d(f,\,g)\, \circ_w\, w_d(g,\,f)$.
The operation $\circ_w$ must be commutative due to the symmetry of $d_w$. From
the remaining properties of Definition \ref{DEF_intervaluation}, property (4)
follows directly from $d(g,\,g) = 0$, while the rest is less obvious.

\begin{Example} \label{EXA_intervaluation_reals}
The standard \imetric on $[0,\,\infty)$ is an intervaluation \imetric with
\beq
w(r,\, s) &:=& 0\,\smax\, (r\,-\,s).
\eeq
However, one may freely choose $\circ_w$ to be addition or maximum. To prove
the cut law for both choices, it suffices to show
\beq
0\,\smax\,(r\,-\,s) &=&
\left(0\,\smax\,\big(r\,-\,(s\,\smax\,t)\big)\right)\;+\;
\left(0\,\smax\,\big((r\,\smin\, t)\,-\, s)\big)\right).
\eeq
For this, we make use of $a\,+\,b\,=\,(a\,\smin\,b)\,+\,(a\,\smax\,b)$ with
$a\,=\,r\,\smin\,s$ and $b\,=\,r\,\smin\,t$, then add $r$ to both sides,
rearrange and apply $x\,-\,(x\,\smin\,y) \,=\, 0\,\smax\,(x\,-\,y)$.
\end{Example}

\begin{Example} \label{EXA_intervaluation_lp}
Let $(X,\,\my)$ be a measure space, $p\,\in\,(1,\infty)$ arbitrary, and
$L$ the lattice of $L^p$-integrable non-negative Lipschitz functions of Lipschitz constant
$\leq 1$. Define
\beq
r\,\circ_w\,s&:=&(r^p\,+\,s^p)^{1/p}, \\
\n{and} \qquad
w(f,\,g) &:=& \sqrt[p]{\int\big|f\,-\,(f\,\smin\, g)\big|^p\,\d\my}\,.
\eeq
As $|r\,-\,(r\smin s)|^p \,+\, |s\,-\,(r\smin s)|^p \,=\, |r\,-\,s|^p$ for all
$r,\,s\,\in\, [0,\,\infty)$, the corresponding
(pseudo-)\imetric is just the $L^p$-\imetric
\beq
d_p(f,g) &=& \sqrt[p]{\int |f\,-\,g|^p\,\d\my}\,.
\eeq
Properties (1)-(3) of Definition \ref{DEF_intervaluation} follow from Example
\ref{EXA_intervaluation_operations}, (4) is trivial. The left cut law can be shown
by pointwise analysis and case distinction ($h\leq g$ vs. $h > g$), the right
cut law follows from Example \ref{EXA_intervaluation_reals} and the Minkowski
inequality. $d_p$ might be a pseudo-metric, depending on $\my$.
\end{Example}


\begin{Example}
Here is a minimal example for a non-intervaluation metric: Take $L=\{a,b,c\}$ with
$a\,<\,b\,<\,c$, and $d(a,\,c) \,=\, 1$, $d(a,\,b)\,=\,2$, $d(b,\,c)\,=\,3$.
Then $w(c,\,a) \,=\, 1$, although $w(c\,\smin\, b,\, a)\,=\,2$ and
$w(c,\,a\,\smax\, b) \,=\, 3$, which both contradict the cut law
and Proposition \ref{PRO_intervaluation}.1, no matter what $\circ_w$ is.
\end{Example}

\begin{Example} \label{EXA_lipschitz_norm}
The Lipschitz constant provides a much more interesting example for a
non-intervaluation metric. Let $X$ be an arbitrary true metric space, and $L$ a
complete lattice of functions $f:\,X\rightarrow \R$ with bounded Lipschitz
constant. The Lipschitz constant of a function $f\in L$ and the corresponding
pseudo-\imetric are given by
\beq
\n{LC}(f) &:=& \sup_{x,\,y\,\in\,X} \frac{\,\big|f(x)\,-\,f(y)\big|\,}{d(x,\,y)} \\
d_\n{LC}\,(f,\,g) &:=& \n{LC}(f\,-\,g).
\eeq
They are used by \cite{Weaver} as ingredient to the utilized norm, called
{\em Lipschitz norm}, which is defined as $||f||_L \,:=\,
||f||_\infty\,\smax\,\n{LC}(f)$. However, neither defines an intervaluation:
Although Weaver shows in his Proposition 1.5.5 that $\n{LC}$ fulfills a modular
inequality for ultravaluations
\beq
\n{LC}(f\,\smax\,g)\;\smax\; \n{LC}(f\,\smin\,g)&\leq& \n{LC}(f)\;\smax\; \n{LC}(g)
\eeq
the inverse inequality is wrong, as there is no bound to $\n{LC}(f)$ by any
combination of $\n{LC}(f\smin g)$ and $\n{LC}(f\smax g)$. To see this, consider
the two-point-space $X\,=\,\{a,b\}$ of diameter $l<1$, and the Lipschitz-functions
$f\,=\,(0,\,l)$ and $g\,=\,(l,\,0)$. Then $\n{LC}(f)\,=\,||f||_L\,=\,1$, but
$\n{LC}(f\smin g)\,=\,\n{LC}(f\smax g)\,=\,0$ and $||\cdot||_L\,=\,l$ in both
cases.

Correspondingly, the cut law is explicitly violated by $d_\n{LC}$, as one can
see when $f$ and $g$ are two different constant functions, and $h$ crosses them
both.
\end{Example}

\goodbreak

We now concentrate on the special case of the supremum metric.

\begin{Proposition} \label{PRO_supremum_intervaluation}
Let $Z$ be a distributive lattice with intervaluation \imetric $d$
(with corresponding $w_d$ and $\circ_d$), with $r\, \circ_d\,s\, =\,r\, \smax\,s$ for all $r,\,s\,\in\,
[0,\,\infty)$. Let $X$ be an arbitrary space, and $L$ a complete lattice of
functions $f:\, X\rightarrow Z$ with pointwise infima and suprema. If
\beq
w_\infty(f,\,g) &:=& \bigsup_{x\in X} w_d\big(f(x),\, g(x)\big)
\eeq
is bounded, it defines an intervaluation \imetric on $L$ with $r\,\circ_\infty\,s\,
=\,r\,\smax\,s$ for all $r,\,s\,\in\, [0,\,\infty)$, which equals the supremum
\imetric $d_\infty$.
\end{Proposition}
\begin{Proof}
The left inequality of the cut law is trivial. For the right side we have to
use that a supremum of sums is less than or equal to a sum of suprema, which in
turn follows from complete distributivity:
\beq
\bigsup_{x\in X} w_d\big(fx,\,gx\big) &\leq& \bigsup_{x\in
X}\left(w_d\big(fx,\, (g\smax h)(x)\big)\;+\;w_d\big((f\smin h)(x),\,
gx\big)\right)\\
&\leq& \bigsup_{x\in X} w_d\big(fx,\,
(g\smax h)(x)\big)\;+\;\bigsup_{x\in X} w_d\big((f\smin h)(x),\, gx\big)
\eeq
\end{Proof}

\begin{Corollary} \label{COR_lipschitz_is_intervaluation}
Let $X$ be any \imetric space. The supremum \imetric $d_\infty$ is an
intervaluation \imetric on the space $\Lip_0 X$ of bounded, non-negative
Lipschitz functions on $X$ with Lipschitz-constant $\leq 1$.
\end{Corollary}
\begin{Proof}
$\Lip X$ is a complete lattice, as one can easily check. We find $r\,\circ_{d_\infty}\, s\;=\; r\,\smax s$ and
\beq
w_{d_\infty}(f,\, g) \quad=\quad \bigsup_{x\in X} \big| f(x) \,-\, (f\,\smin\, g)(x)\big|
\quad=\quad 0\,\smax\, \bigsup_{x\in X} \big(f(x) \,-\, g(x)\big),
\eeq
which is the intervaluation \imetric of Proposition
\ref{PRO_supremum_intervaluation} applied to Example
\ref{EXA_intervaluation_reals}.
\end{Proof}

\subsection{Topological Aspects}

We finally take a look at the subset $\mli(L)$ of all completely
$d$-irreducible elements of a complete lattice $L$ with
intervaluation metric $d$.

\begin{Proposition} \label{PRO_ml-irreducible_complete}
Let $L$ be a complete lattice with intervaluation \imetric $d$, and let $L$ be
metrically complete. Then $\mli(L)$ is topologically closed.
\end{Proposition}
\begin{Proof}
Let $(p_n)\subset \mli(L)$, $n\in \N^*$ be some sequence of completely $d$-irre\-ducible elements
converging to $p\,\in\, L$, and $(f_j)_{j\in J}$ any non-empty family in $L$. Then
for any $n\,\in\, \N^*$ holds
\beq
d\,\left(p,\,\bigsup f_j\right) &\geq& d\,\left(p_n,\,\bigsup f_j\right) \;-\;
d(p,\, p_n)\\
&\geq& \biginf d(p_n,\, f_j) \;-\; d(p,\,p_n)\\
&\geq& \biginf \big(d(p,\, f_j) \,-\, d(p,\, p_n)\big) \;-\; d(p,\,p_n)\\
&\geq& \biginf d(p,\, f_j) \;-\; \underbrace{2~d(p,\,p_n)}_{\rightarrow\; 0},
\eeq
i.e.\ the element $p$ is completely $d$-irreducible.

\end{Proof}

\begin{Definition} \label{DEF_order_base}
Let $L$ be a lattice with \imetric $d$, $R\geq 0$ arbitrary. We define an {\em
$R$-base} of $L$ to be a subset $B\,\subset\, L$ such that for any $f\in L$
there is $(b_j)_{j\in J}\,\subset\, B$, $J$ an arbitrary non-empty index set,
such that $d(f,\,\bigsup_{j\in J}\,b_j) \,\leq\, R$. A {\em base} simply is a
$0$-base.
\end{Definition}


\begin{Proposition} \label{PRO_irred_in_base}
Consider an $R$-base $B$ of a complete lattice $L$ with intervaluation \imetric
$d$, $R\geq 0$. Then for each $\delta\,>\,0$, $\mli(L)$ is in the
$(R\,+\,\delta)$-ball around $B$. In particular, if $R\,=\,0$, $\mli(L)$ lies in
the metrical closure of $B$.
\end{Proposition}
\begin{Proof}
Let $p\,\in\, \mli(L)$ be arbitrary. As $B$ is an $R$-base, there are $b_j\,\in\,
B$, $j\,\in\, J\,\neq\,\emptyset$, such that
\beq
d\left(p,\, \bigsup_{j\in J} b_j\right) &\leq& R.
\eeq
From Definition \ref{DEF_ml-irreducible} we infer that there is a sequence
$(c_k)\,\subset\, B$, $k\,\in\, K\,\subset\, J$ whose distances to $p$
converge to $R$. If $R\,=\,0$, the sequence $(c_j)$ metrically converges to
$p$.
\end{Proof}

Propositions \ref{PRO_ml-irreducible_complete} and \ref{PRO_irred_in_base}
might help in identifying all completely $d$-irreducible elements of a
concretely given lattice. 

\begin{Example}
It is easy to see that, if $B$ is a base, and $b\,\in\, B$ not a join-irreducible
element, then $B\setminus\{b\}$ is a base as well (if $b\,=\,f\,\smax\, g$, $f$
and $g$ are joins of elements of $B$, and as $f,\,g\,<\,b$, $b$ is not part of
these joins). Using the Lemma of Zorn, it is possible to deduce that the subset of
all join-irreducible elements constitutes a base for any sufficiently nice
lattice.

Unfortunately, this is not the case with $d$-irreducible elements: Let $L'$ be
the completely distributive complete lattice $[0,\,3]\times[0,\,2]$ with
componentwise supremum and infimum, and with supremum metric. Then consider the
sublattice $L\subset L'$ formed by the five elements
\beq
L&:=&\{(0,\,0),\; (1,\,0),\; (0,\,1),\; (1,\,1),\; (2,\,2)\}.
\eeq
We find $\mli(L)\,=\,\{(0,\,0),\, (1,\,0),\, (0,\,1)\}$, as
$(1,\,1)\,=\, (1,\,0)\,\smax\,(0,\,1)$. $p\,=\,(2,\,2)$ is join-irreducible in
this lattice, but not $d$-irreducible: Take $f_1\,=\,(1,\,0)$, $f_2\,=\,(0,\,1)$,
then $\biginf d(p,\,f_j)\, =\, 2$, but $d(p,\,\bigsup f_j)\, =\, 1$.
Nevertheless, $(2,\,2)$ must be part of any 0-base of $L$.
\end{Example}

Georg-August-Universit\"at G\"ottingen, Germany\newline
eMail \verb|lochmann@uni-math.gwdg.de|


\begin{thebibliography}{9999}

\bibitem[Bi1]{Birkhoff}G. Birkhoff, {\it Lattice Theory}, American Mathematical
  Society Colloquium Publications Vol. XXV, 2nd ed. (1948) and 3rd ed. (1960)

\bibitem[Bi2]{Birkhoff_on_von_Neumann}G. Birkhoff, {\it Von Neumann and Lattice
  Theory}, Bull. Amer. Math. Soc. {\bf 64}, Nr 3, Part 2 (1958) 50--56, \\
  \href{http://www.ams.org/bull/1958-64-03/S0002-9904-1958-10192-5/S0002-9904-1958-10192-5.pdf}
  {http://www.ams.org/bull/1958-64-03/S0002-9904-1958-10192-5/} \\  
  \mbox{}\qquad\quad \href{http://www.ams.org/bull/1958-64-03/S0002-9904-1958-10192-5/S0002-9904-1958-10192-5.pdf}
  {S0002-9904-1958-10192-5.pdf}

\bibitem[dH]{dlHarpe}P. de la Harpe, {\it Topics in Geometric Group Theory}, The
  University of Chicago Press (2000)

\bibitem[Gl]{Glivenko}V. Glivenko, {\it G\'eometrie des syst\`emes de chosen
  norm\'ees}, Am. Jour. of Math. {\bf 58} (1936) 799--828

\bibitem[Lo1]{Lochmann_dissertation}A. Lochmann, {\it Rough Isometries of Order
  Lattices and Groups}, Nieders\"achsische Staats- und Universit\"atsbibliothek,
  Doctoral Thesis,
  \href{http://webdoc.sub.gwdg.de/diss/2009/lochmann/}{http://webdoc.sub.gwdg.de/diss/2009/lochmann/}

\bibitem[Lo2]{Lochmann_rough_isometries_of_lipschitz_function_spaces}
  A. Lochmann, {\it Rough Isometries of Lipschitz Function Spaces},
  preprint at \href{http://arxiv.org/abs/0710.1109}{http://arxiv.org/abs/0710.1109}

\bibitem[Mn]{Monjardet}B. Monjardet, {\it Metrics on partially ordered sets ---
  a survey}, Discrete Mathematics {\bf 35} (1981) 173--184

\bibitem[Ro]{Rockafellar}R. T. Rockafellar, {\it Convex Analysis}, Princeton
  University Press (1970)

\bibitem[vN]{von_Neumann_continuous_geometries}J. von Neumann, {\it Lectures on
  continuous geometries}, Princeton 1936-1937 (2 vols.), in particular chapter
  XVII

\bibitem[Wr]{Werner}D. Werner, {\it Funktionalanalysis}, Springer (2005)

\bibitem[Wv]{Weaver}N. Weaver, {\it Lipschitz Algebras}, World Scientific (1999)

\end{thebibliography}
\end{document}